# THE PROBLEM OF FINDING THE GLOBAL MINIMUM OF A POLYNOMIAL ON A COORDINATE PARALLELEPIPED

**V. N. Nefedov**

**Abstract.** The paper addresses the problem of finding the global minimum of a polynomial on a coordinate parallelepiped or on certain other sets with a simple structure. To solve the problem, the Variable Elimination Algorithm (VEA) for polynomial optimization problems is employed, as described in [1, 2]. It should be noted that the proposed method can be of practical significance only for problems of low dimensionality and for polynomials of not very high degree. This is due to the fact that, when implementing the VEA, a large number of alternative problems may arise concerning the determination of real roots of polynomials, and these polynomials may have very high degrees. Although, in theory, the VEA is applicable to any polynomial optimization problems, there are often practically relevant problems that can be successfully solved using the proposed method.

**Introduction**

When solving certain applied problems, we arrive at the problem of finding the global minimum of a polynomial $p(\mathbf{x})$, where $\mathbf{x} = (x_1, \dots, x_m) \in \mathbb{R}^m$, on a coordinate parallelepiped $\Pi$ or on a set $X \subset \mathbb{R}^m$ with a simple structure (for example, described by a finite collection of inequalities using linear or quadratic functions). If the polynomial $p(\mathbf{x})$ is not convex (and the set $X$ may also turn out to be non-convex), we face a computationally difficult problem. For instance, applying traditional gradient methods will only allow us to find some point of local minimum.

In the case of a coordinate parallelepiped, it may generally be necessary to compute the values of $p(\mathbf{x})$ at the nodes of a uniform grid $\Pi^\tau \subset \Pi$ with a small step $\tau > 0$ along each coordinate. If $X$ differs from a coordinate parallelepiped, constructing a similar grid $X^\tau \subset X$ is yet another challenging task. Moreover, to assess how accurately we have determined the minimum value by computing $\min\{p(\mathbf{x}) \mid \mathbf{x} \in \Pi^\tau\}$, we need to know the Lipschitz constant — that is, a value $L > 0$ such that

$$\forall \mathbf{x}, \tilde{\mathbf{x}} \in \Pi \quad | p(\mathbf{x}) - p(\tilde{\mathbf{x}}) | \le L \parallel \mathbf{x} - \tilde{\mathbf{x}} \parallel, \qquad (0.1)$$

where $\parallel \mathbf{x} \parallel = \max\{| x_i | : i = 1, \dots, m\}$. However, to obtain at least an approximate value of this constant, we will need to solve another problem:

$$\sum_{i=1}^{m} \left| \frac{\partial p(\mathbf{x})}{\partial x_i} \right| \to \max\,(= \bar{L});\, \mathbf{x} \in \Pi.$$

Any number greater than or equal to $\bar{L}$ can be taken as the Lipschitz constant $L$ satisfying (0.1). Thus, we arrive at yet another global optimization problem. At the same time, the value of $\bar{L}$ may turn out to be very large, making its practical application infeasible.

The aforementioned difficulties can sometimes be circumvented by using the Variable Elimination Algorithm for polynomial optimization problems (see [1, 2]). This algorithm (hereinafter referred to as VEA) is applied to a polynomial problem in canonical form, i.e., to a problem of the form

$$x_1 \to \min ; p_i(\mathbf{x}) = 0,\ i = 1, \dots, l, \qquad (0.2)$$

where $p_1(\mathbf{x}), \dots, p_l(\mathbf{x})$ are polynomials. The algorithm consists in sequentially eliminating the variables $x_m, x_{m-1}, \dots, x_2$ until we obtain a finite collection of alternative problems of the form

$$x_1 \to \min ; h(x_1) = 0, \qquad (0.3)$$

where $h(x_1)$ is a univariate polynomial in $x_1$. In this case, any local minimum value (including the global one) in problem (0.2) will be a local minimum value in one of the alternative problems of type (0.3), and thus will be a real root of the corresponding polynomial $h(x_1)$. Consequently, we can obtain a finite set $U$ of numbers (the real roots of the polynomials from the entire collection of alternative problems of type (0.3)), which includes all local minimum values of problem (0.2). The local minimum values (including the global one) are computed **with the same accuracy** as the accuracy with which we can find the real roots of univariate polynomials.

Next, by iterating over the found values $x_1 \in U$, we can reconstruct, via the chain of alternative problems obtained as a result of the sequential variable elimination, all components of any point of local (and hence global) minimum. This ensures good accuracy not only in terms of the value of the minimized function but also in terms of the distance to the set of global minimum points Argmin $p(\Pi)$.

Furthermore, as will become clear later, we will obtain a solution to a two-stage lexicographic optimization problem: specifically, we will find a point $\mathbf{x}^* \in$ Argmin $p(\Pi)$ with the minimal value of the $x_1$ component. Similarly, we can find a point from Argmin $p(\Pi)$ with the maximal value of this component, as well as with the minimal and maximal values of all other components. By proceeding in this manner, we can construct the smallest coordinate parallelepiped containing the set Argmin $p(\Pi)$, or make sure that the global minimum is attained at a single point.

The main result of this paper consists of statements (see Lemma 1, Theorem 1) that allow reducing the original problem

$$p(\mathbf{x}) \to \min ; \mathbf{x} \in \Pi \qquad (0.4)$$

to a collection of similar problems on the faces of the parallelepiped $\Pi$ (thereby reducing the problem dimension by 1), as well as to a problem (in canonical form) of the form

$$x_1 \to \min ; \frac{\partial p(\mathbf{x})}{\partial x_i} = 0,\ i = 1, \dots, m, \qquad (0.5)$$

to which the Variable Elimination Algorithm (VEA) is fully applicable. Here, problem (0.5) corresponds to the case where all global minimum points of problem (0.4) are interior points of the parallelepiped $\Pi$.

One could apply the VEA to problem (0.4) without relying on these results. However, in that case, the scenario where all global minimum points of problem (0.4) are interior points of $\Pi$ would correspond to a canonical-form problem of the type

$$u \to \min\,;\, u - p(\mathbf{x}) = 0. \qquad (0.6)$$

Yet, unlike problem (0.5), problem (0.6) introduces an additional variable, which typically leads to a significant increase in the complexity of solving the problem (due to the growth in the number of alternative problems and in the degrees of the polynomials involved in them).

Note that the VEA can be applied to solve certain problems of type (0.4) in which the objective function is not a polynomial. For example, consider the problem

$$f(\mathbf{x}) \to \min\,;\, \mathbf{x} \in \Pi, \qquad (0.7)$$

where $f(\mathbf{x}) = \frac{p_1(\mathbf{x})}{p_2(\mathbf{x})}$ is a rational function, $p_1(\mathbf{x}), p_2(\mathbf{x})$ are polynomials, and $p_2(\mathbf{x}) \neq 0$ for all $\mathbf{x} \in \Pi$. Then, this problem can be transformed into an equivalent one:

$$u \to \min\,;\, up_2(\mathbf{x}) - p_1(\mathbf{x}) = 0,\ \mathbf{x} \in \Pi. \qquad (0.8)$$

When solving problem (0.8), the case where all global minimum points of problem (0.7) are interior points of the parallelepiped $\Pi$ corresponds to a polynomial problem in canonical form:

$$u \to \min\,;\, up_2(\mathbf{x}) - p_1(\mathbf{x}) = 0, \qquad (0.9)$$

to which the VEA is fully applicable.

Let us also consider another case. Suppose we are dealing with a problem of the form

$$|\, p_1(\mathbf{x}) \,| + |\, p_2(\mathbf{x}) \,| \to \max\,;\, \mathbf{x} \in \Pi, \qquad (0.10)$$

where $p_1(\mathbf{x}), p_2(\mathbf{x})$ are polynomials. As we have already seen, problems of this kind arise when computing the Lipschitz constant of a polynomial. Then, problem (0.10) can be reduced to considering four problems:

$$(-1)^i p_1(\mathbf{x}) + (-1)^j p_2(\mathbf{x}) \to \min\,;\, \mathbf{x} \in \Pi, \qquad (0.11)$$

where $i, j \in \{0,1\}$. The solution to problem (0.10) will be one of the solutions to problems (0.11) that yields the maximum value of the objective function $|\, p_1(\mathbf{x}) \,| + |\, p_2(\mathbf{x}) \,|$. Problems (0.11) are problems of type (0.4).

The problem of finding the global minimum of a polynomial on a coordinate parallelepiped belongs to the class of computationally difficult global optimization problems, especially in the case of non-convex polynomials and non-convex feasible sets. Traditional gradient methods allow finding only local minima, while the use of grid-based approaches entails the need to estimate the Lipschitz constant and is associated with an exponential growth in the number of grid nodes as the problem dimension increases — which makes such approaches poorly suited even for problems of moderate dimensionality.

Alternative modern approaches have their own strengths and weaknesses. For instance, the Lasserre hierarchy, based on sums of squares (SOS) and moment matrices, provides theoretically grounded lower bounds for the objective functional and, in the limit, yields the exact value of the global minimum. However, at finite levels of the hierarchy, it delivers only approximate estimates, and the computational complexity grows rapidly with an increase in the number of variables and the degree of the polynomial [3]. Further development and systematization of SOS optimization ideas — including the connection with moment matrices and practical aspects of application — are described in detail in the survey [4].

Branch-and-bound methods and interval optimization are capable of guaranteeing the global optimality of the solution with a controlled error; however, in the worst case, they also face exponential complexity. Their effectiveness largely depends on the quality of the cuts and the heuristics used for partitioning the domain (see [5–7], as well as the survey [8]).

In this context, the Variable Elimination Algorithm (VEA), considered in this paper, offers a different approach: by sequentially eliminating variables, the problem is reduced to a finite set of univariate polynomial equations, the roots of which determine all local (including the global) minimum points. An important advantage of the method is the explicit consideration of the domain boundary through the recursive examination of the faces of the parallelepiped with a reduction in dimensionality. Moreover, the approach can be extended to rational functions and problems involving absolute values by transforming them into equivalent polynomial formulations.

At the same time, the practical applicability of the VEA is limited to problems of small dimensionality and moderate polynomial degree due to the rapid growth in the number of alternative subproblems and the degrees of the intermediate polynomials. Thus, the method occupies a specific niche among tools for global polynomial optimization: it provides a fairly precise algebraic description of all candidate points for the minimum in problems with a simple structure, but it is inferior in scalability to a number of modern numerical approaches — in particular, to the SOS methods developed in [3, 4].

## 1. The Variable Elimination Algorithm in Polynomial Optimization Problems (VEA)

We provide a brief description of the VEA. A detailed exposition can be found in [1, 2]. Theoretically, this algorithm is applicable to any polynomial problems — that is, to problems with an arbitrary finite number of variables and any number of polynomial constraints of equality or inequality type. The result of applying this algorithm is a finite set containing all values of the local extrema of the original problem (this set may also include some other numbers). This set is the collection of real roots of the univariate polynomials obtained as a result of the VEA operation.

It should be noted immediately that, since the implementation of the VEA may give rise to a large number of alternative problems for finding the real roots of polynomials — and these polynomials may have very high degrees — the VEA can be of practical significance only for problems of low dimensionality and for polynomials in the initial conditions of not very high degree.

A polynomial programming problem is a problem of the form

$$p_0(x_1, \dots, x_m) \to \min \ \text{(or max)}; \qquad (1.1)$$

$$p_i(x_1, \dots, x_m) \le 0,\ i = 1,2, \dots, l_1; p_j(x_1, \dots, x_m) = 0,\ j = l_1 + 1, \dots, l, \quad (1.2)$$

where $(x_1, \dots, x_m) \in \mathbb{R}^m$ and $p_i, p_j$ are polynomials.

We will consider problem (1.1), (1.2) in canonical form:

$$x_1 \to \text{extr (i.e., either min or max )}; \quad (1.3)$$

$$p_i(x_1, \dots, x_m) = 0,\ i = 1,2, \dots, l, \quad (1.4)$$

to which the VEA will be applied. Problem (1.1), (1.2) can be easily reduced to a problem of type (1.3), (1.4) by introducing additional variables and constraints [1,2].

To make the reasoning more illustrative, let us consider the case where there are three variables and the problem takes the form:

$$x \to \text{extr}; p_i(x, y, z) = 0,\ i = 1,2, \dots, l. \quad (1.5)$$

We describe the algorithm for transforming problem (1.5) into several alternative problems of the same type, such that the set of local extremum values in problem (1.5) is contained in the set of local extremum values of the alternative problems, and these alternative problems do not involve the variable $z$.

By the degree of the variable $z$ in problem (1.5), we mean the maximum degree of this variable in the polynomials $p_i(x, y, z), i = 1,2, \dots, l$.

We will consider the following three cases.

**Case 1.** Let $l = 1$, i.e., problem (1.5) takes the form

$$x \to \text{extr}; p_1(x, y, z) = 0, \quad (1.6)$$

and the degree of the variable $z$ in the polynomial $p_1(x, y, z)$ is $n \ge 1$. Then, for any point $(x^*, y^*, z^*) \in \mathbb{R}^3$ that is a local extremum point in problem (1.6), the following holds (see [2, Lemma 1]):

$$\frac{\partial p_1(x^*, y^*, z^*)}{\partial z} = 0,$$

and consequently, $(x^*, y^*, z^*)$ is a local extremum point in the problem

$$x \to \text{extr}; p_1(x, y, z) = 0,\ \frac{\partial p_1(x, y, z)}{\partial z} = 0,$$

which is equivalent to the problem

$$x \to \text{extr}; n p_1(x, y, z) - z \frac{\partial p_1(x, y, z)}{\partial z} = 0,\ \frac{\partial p_1(x, y, z)}{\partial z} = 0 \quad (1.7)$$

(these problems have identical feasible sets of possible solutions). At the same time, the degree of the variable $z$ in problem (1.7) is one unit less than in the original problem (1.6).

**Case 2.** Let $l > 1$, and the variable $z$ appears in only one constraint, i.e., problem (1.5) has the form

$$x \to \text{extr}; p_1(x, y, z) = 0,\ p_i(x, y) = 0,\ i = 2, \dots, l. \quad (1.8)$$

Furthermore, let the degree of the variable $z$ in the polynomial $p_1(x, y, z)$ be $n \ge 1$. In this case, for any point $(x^*, y^*, z^*) \in \mathbb{R}^3$ that is a local extremum point in problem (1.8), one of the following holds:

- **(a)** either $\frac{\partial p_1(x^*, y^*, z^*)}{\partial z} = 0$, and then $(x^*, y^*, z^*)$ is a local extremum point in the problem

$$x \to \text{extr}; p_1(x, y, z) = 0,\ \frac{\partial p_1(x, y, z)}{\partial z} = 0,\ p_i(x, y) = 0,\ i = 2, \dots, l,$$

which is equivalent to the problem

$$x \to \text{extr}; np_1(x,y,z) - z\frac{\partial p_1(x,y,z)}{\partial z} = 0, \frac{\partial p_1(x,y,z)}{\partial z} = 0, \quad (1.9)$$
$$p_i(x,y) = 0, i = 2, \ldots, l,$$

where the degree of the variable $z$ is one unit less than in problem (1.8);

- **(b)** or $\frac{\partial p_1(x^*,y^*,z^*)}{\partial z} \neq 0$, and then $(x^*, y^*)$ is a local extremum point in the problem

$$x \to \text{extr}; p_i(x,y) = 0, \ i = 2, \ldots, l, \quad (1.10)$$

in which the variable $z$ has been eliminated.

Thus, in Case 2, we move from the single problem (1.8) to two alternative problems (1.9) and (1.10) such that any local extremum value in problem (1.8) is a local extremum value in at least one of the problems (1.9) or (1.10), and the degree of the variable $z$ in these problems is either one unit less than in problem (1.8) or equal to 0 (i.e., the variable $z$ is eliminated).

**Case 3.** To consider Case 3, we need a special operation on polynomials.

We introduce the operation $R_1^z$ — the elementary reduction of the degree of the variable $z$ in the polynomial

$$p(x,y,z) = p_{n_1}(x,y)z^{n_1} + p_{n_1-1}(x,y)z^{n_1-1} + \cdots + p_1(x,y)z + p_0(x,y)$$

by means of the polynomial

$$q(x,y,z) = q_{n_2}(x,y)z^{n_2} + q_{n_2-1}(x,y)z^{n_2-1} + \cdots + q_1(x,y)z + q_0(x,y),$$

where

$$p_{n_1}(x,y), p_{n_1-1}(x,y), \ldots, p_1(x,y), p_0(x,y),$$
$$q_{n_2}(x,y), q_{n_2-1}(x,y), \ldots, q_1(x,y), q_0(x,y)$$

are polynomials, $p_{n_1}(x,y) \not\equiv 0$, $q_{n_2}(x,y) \not\equiv 0$, and $n_1 \geq n_2 \geq 1$. The result of this operation is the polynomial:

$$R_1^z(p,q) = q_{n_2}(x,y)p(x,y,z) - z^{n_1-n_2}p_{n_1}(x,y)q(x,y,z).$$

Let $\tilde{p}(x,y,z) = R_1^z(p,q)$. Then

$$\tilde{p}(x,y,z) = \tilde{p}_{n_1-1}(x,y)z^{n_1-1} + \cdots + \tilde{p}_1(x,y)z + \tilde{p}_0(x,y),$$

i.e., the degree of the variable $z$ in the polynomial $\tilde{p}(x,y,z)$ does not exceed $n_1 - 1$.

In the case $n_1 \geq n_2 + 1$, the polynomials $p(x,y,z)$ and $q(x,y,z)$ can also be subjected to the operation $R_2^z$ — a double reduction of the degree of the variable $z$ in the polynomial $p(x,y,z)$ by means of the polynomial $q(x,y,z)$, defined by the equality

$$R_2^z(p,q) = R_1^z(R_1^z(p,q),q).$$

Note that the degree of the variable $z$ in the polynomial $\check{p}(x,y,z) = R_2^z(p,q)$ does not exceed $n_1 - 2$. Similarly, we can introduce operations for triple, …, $(n_1 - n_2 + 1)$-fold elementary reduction of the degree of the variable $z$ in the polynomial $p(x,y,z)$ by means of the polynomial $q(x,y,z)$; we will denote these operations by

$$R_1^z, R_2^z, \ldots, R_{n_1-n_2+1}^z,$$

respectively. In this case, the degree of the variable $z$ in the polynomial $R_j^z(p,q)$, where $1 \le j \le n_1 - n_2 + 1$, does not exceed $n_1 - j$.

We generalize the introduced operation to the case where $0 \le n_1 < n_2$ (i.e., it is even permissible for the polynomial $p(x,y,z)$ not to depend on the variable $z$ at all); in this case, we set $R_1^z(p,q) = p(x,y,z)$, and consequently, for this case, $\forall k \in \mathbb{N}$, $R_k^z(p,q) = p(x,y,z)$.

**Remark 1.** In what follows, the operation $R_1^z$ will be applied repeatedly to eliminate the variable $z$ from a sequence of alternative problems. At the same time, the degrees of the remaining variables may grow rapidly. To reduce the rate of this growth as much as possible, we can further modify the operation $R_1^z$ as follows. If

$$g(x,y) = \mathrm{GCD}(p_{n_1}(x,y), q_{n_2}(x,y)) \not\equiv 1, p_{n_1}(x,y) = \bar{p}_{n_1}(x,y)g(x,y),$$
$$q_{n_2}(x,y) = \bar{q}_{n_2}(x,y)g(x,y),$$

then we set

$$R_1^z(p,q) = \bar{q}_{n_2}(x,y)p(x,y,z) - z^{n_1-n_2}\bar{p}_{n_1}(x,y)q(x,y,z).$$

Now we can return to the description of Case 3.

Let, in problem (1.5), $l \ge 2$, and $k \ge 2$ polynomials among $p_i(x,y,z)$, $i = 1,2,\dots,l$, have a positive degree with respect to the variable $z$, and the maximum degree of the variable $z$ in these polynomials is $n \ge 1$. Moreover, for definiteness, let the variable $z$ enter the polynomial $p_1(x,y,z)$ with the smallest positive degree $n_1 \ge 1$ (see Remark 2 regarding the possibility of uniquely selecting the polynomial $p_1(x,y,z)$), where $1 \le n_1 \le n$, i.e.

$$p_1(x,y,z) = p_{1,n_1}(x,y)z^{n_1} + p_{1,n_1-1}(x,y)z^{n_1-1} + \dots + p_{1,1}(x,y)z + p_{1,0}(x,y),$$

with $p_{1,n_1}(x,y) \not\equiv 0$.

Let $(x^*, y^*, z^*)$ be an arbitrary local extremum point in problem (1.5). Then one of the following holds:

- **(a)** either $p_{1,n_1}(x^*,y^*) \neq 0$, and then $(x^*,y^*,z^*)$ is a local extremum point in the problem

$$\begin{gathered} x \to \mathrm{extr}; p_1(x,y,z) = 0, \\ \tilde{p}_i(x,y,z) = R_{n-n_1+1}^z(p_i,p_1) = 0,\ i = 2,\dots,l, \end{gathered} \tag{1.11}$$

(here, the operation $R_{n-n_1+1}^z$ is applied in the generalized sense noted in Remark 1: it leaves unchanged polynomials that do not depend on $z$ or have a degree in $z$ less than $n_1$, while transforming all polynomials with degree greater than or equal to $n_1$ into polynomials with degree less than or equal to $n_1 - 1$);

- **(b)** or $p_{1,n_1}(x^*,y^*) = 0$, and then $(x^*,y^*,z^*)$ is a local extremum point in the problem

$$\begin{gathered} x \to \mathrm{extr}; p_{1,n_1}(x,y) = 0, \\ \tilde{p}_1(x,y,z) = p_{1,n_1-1}(x,y)z^{n_1-1} + \dots + p_{1,1}(x,y)z + p_{1,0}(x,y) = 0, \\ p_i(x,y,z) = 0,\ i = 2,\dots,l. \end{gathered} \tag{1.12}$$

Thus, in each of Cases 1–3, we transform the original problem (1.5) into new alternative problems of the same type (1.5) (one problem in Case 1, two in Case 2, and two in Case 3). From these, in turn, we can proceed to further alternative problems, etc. In other words, we can organize a multi-stage branching process of the original problem (1.5) into a multitude of alternative problems of the same form (i.e., again of type (1.5)), with a successive reduction of the degree of the variable $z$ until it is completely eliminated.

At each stage of this process, we have a finite collection of alternative problems such that the union of the sets of local extremum values in these problems contains the set of local extremum values of the original problem (1.5). The goal of this process is to obtain a collection of alternative problems with the variable $z$ eliminated. A similar process will then allow us to eliminate the variable $y$, yielding a collection of alternative problems of type (1.5) with a single variable $x$, i.e., problems of the form

$$x \to \text{extr}; p_i(x) = 0,\ i = 1, \dots, l, \qquad (1.13)$$

where $p_i(x)$ are univariate polynomials, $i = 2, \dots, l$. Further, in Remark 2 it will be shown that for each of the alternative problems of the form (1.13), i.e., with the single variable $x$, the condition $l = 1$ will hold.

**Remark 2.** Note that, after obtaining an alternative problem of type (1.13), it can be simplified. First, the solution set of such a problem is the set of common real roots of the univariateй polynomials $p_i(x)$, $i = 1, \dots, l$, which equals the set of real roots of the polynomial $p(x) = \text{GCD}\,(p_1(x), \dots, p_l(x))$. As a result of this simplification, only a single polynomial $p(x)$ remains in a problem of type (1.13). To further reduce the degree of the polynomial $p(x)$, we can additionally reduce the multiplicities of the real roots to 1 by considering, instead of $p(x)$, the polynomial

$$\tilde{p}(x) = \frac{p(x)}{\text{GCD}\left(p(x), \frac{dp(x)}{dx}\right)}. \qquad (1.14)$$

Similar simplifications can also be applied to alternative problems with a larger number of variables. For example, suppose we have an alternative problem involving two variables $x$ and $y$, i.e., a problem of the form

$$x \to \text{extr}; p_i(x, y) = 0, i = 1, \dots, l. \qquad (1.15)$$

Let us expand the polynomials $p_i(x, y)$, $i = 1, \dots, l$, in powers of the variable $y$:

$$p_i(x, y) = a_{i,s_i}(x) y^{s_i} + \dots + a_{i,1}(x) y + a_{i,0}(x),$$

where $s_i \geq 0$ and $a_{i,s_i}(x) \not\equiv 0$. Let us set

$$q(x) = \text{GCD}\,(a_{1,s_1}(x), \dots, a_{1,0}(x), \dots, a_{l,s_l}(x), \dots, a_{l,0}(x)).$$

If $q(x) \not\equiv 1$, then, if $(x^*, y^*)$ is a local extremum point in problem (1.15), either $q(x^*) = 0$, and then $x^*$ is a local extremum point in the problem

$$x \to \text{extr}; q(x) = 0, \qquad (1.16)$$

or $q(x^*) \neq 0$, and then $(x^*, y^*)$ is a local extremum point in the problem

$$x \to \text{extr}; \tilde{p}_i(x, y) = 0, i = 1, \dots, l, \qquad (1.17)$$

where $\tilde{p}_i(x, y) = p_i(x, y)/q(x)$, $i = 1, \dots, l$. Thus, in the case $q(x) \not\equiv 1$, we reduce problem (1.15) to the simpler alternative problems (1.16) and (1.17). Moreover, further simplification of the polynomial $q(x)$ is possible according to formula (1.14).

**Remark 3.** If, in Case 3, there are several polynomials with the smallest degree of the variable $z$ equal to $n_1$, we choose among them the polynomial with the smallest degree of the variable $y$ in the polynomial $p_{1,n_1}(x,y)$ (in order to minimize the maximum degree of the variable $y$ in the alternative problems obtained as a result of eliminating the variable $z$). If there are still several such polynomials, we choose among them the one with the smallest degree of the variable $x$ in the polynomial $p_{1,n_1}(x,y)$.

## 2. Application of the Variable Elimination Algorithm in Polynomial Optimization Problems for Finding the Global Minimum of a Polynomial on a Coordinate Parallelepiped

We consider the problem (of finding the global minimum):

$$p(x,y,z) \to \min;\\ (x,y,z) \in \Pi = \{(x,y,z) \in \mathbb{R}^3 \mid a_1 \le x \le a_2, b_1 \le y \le b_2, c_1 \le z \le c_2\}, \quad (2.1)$$

where $p(x,y,z)$ is a polynomial, and $a_1 < a_2$, $b_1 < b_2$, $c_1 < c_2$.

Alongside problem (2.1), we will consider an auxiliary problem:

$$x \to \min ; \frac{\partial p(x,y,z)}{\partial x} = 0,\ \frac{\partial p(x,y,z)}{\partial y} = 0,\ \frac{\partial p(x,y,z)}{\partial z} = 0. \quad (2.2)$$

We will use the following statement, whose proof is given in Section 3.

**Lemma 1.** Let $(x^*,y^*,z^*)$ be a local minimum point of the polynomial $p(x,y,z)$, and at the same time let $(x^*,y^*,z^*)$ be a local minimum point in the problem:

$$x \to \min ; p(x,y,z) = p(x^*,y^*,z^*). \quad (2.3)$$

Then $(x^*,y^*,z^*)$ is a local minimum point in problem (2.2).

In addition, we will need the following simple

**Statement 1.** Let $U, X \subseteq \mathbb{R}^3$, let $f(x,y,z)$ be a function defined on $\mathbb{R}^3$, with $U, X \neq \emptyset$, $U \subset \operatorname{int} X$, and suppose that

$$\forall (x,y,z) \in X \setminus U \quad f(x,y,z) > f(x^*,y^*,z^*). \quad (2.4)$$

Moreover, let $(x^*,y^*,z^*)$ also be a local minimum point in the problem

$$x \to \min ; (x,y,z) \in U.$$

Then $(x^*,y^*,z^*)$ is a local minimum point in the problem

$$x \to \min ; f(x,y,z) = f(x^*,y^*,z^*).$$

**Proof.** Suppose the contrary. Then there exists a sequence $(x^{(n)},y^{(n)},z^{(n)}) \in \mathbb{R}^3$, $n = 1,2,...$, such that

$$\left(x^{(n)},y^{(n)},z^{(n)}\right) \to (x^*,y^*,z^*) \text{ as } n \to \infty,$$

$$f(x^{(n)},y^{(n)},z^{(n)}) = f(x^*,y^*,z^*), x^{(n)} < x^*, n = 1,2,\dots . \quad (2.5)$$

Since $(x^*,y^*,z^*) \in U$ and $U \subset \operatorname{int} X$, for simplicity of notation we assume that

$$(x^{(n)},y^{(n)},z^{(n)}) \in X, n = 1,2,...,$$

and hence (see the conditions of the statement), in view of (2.4) and (2.5), we have:

$$(x^{(n)},y^{(n)},z^{(n)}) \in U, n = 1,2,...$$

But then, using (2.5) again, we arrive at a contradiction with the fact that $(x^*, y^*, z^*)$ is a local minimum point in the problem $x \to \min\,;\,(x, y, z) \in U$. ∎

We denote

$$\text{int}\,\Pi = \{(x, y, z) \in \mathbb{R}^3 \mid a_1 < x < a_2, b_1 < y < b_2, c_1 < z < c_2\}$$

as the set of interior points of the parallelepiped $\Pi$,

$$\begin{aligned}\text{bd}\,\Pi = \{(x, y, z) \in \Pi \mid x = a_1\} \cup \{(x, y, z) \in \Pi \mid x = a_2\}\\ \cup \{(x, y, z) \in \Pi \mid y = b_1\} \cup \{(x, y, z) \in \Pi \mid y = b_2\}\\ \cup \{(x, y, z) \in \Pi \mid z = c_1\} \cup \{(x, y, z) \in \Pi \mid z = c_2\}\end{aligned}$$

as the set of boundary points of the parallelepiped $\Pi$,

$$\text{Argmin}\,p(\Pi) = \{(x, y, z) \in \Pi \mid p(x, y, z) = p^* = \min\,p(\Pi)\}$$

as the set of global minimum points of the polynomial $p(x, y, z)$ on the parallelepiped $\Pi$, and

$$\begin{aligned}x^* &= \min\{x \mid (x, y, z) \in \Pi, p(x, y, z) = p^* = \min\,p(\Pi)\}\\ &= \min\{x \mid (x, y, z) \in \text{Argmin}\,p(\Pi)\}.\end{aligned} \tag{2.6}$$

**Statement 2.** Two cases are possible:

1. $\text{bd}\,\Pi \cap \text{Argmin}\,p(\Pi) \neq \emptyset$;
2. $\text{bd}\,\Pi \cap \text{Argmin}\,p(\Pi) = \emptyset$, and the value $x^*$ is one of the local minimum values in problem (2.2).

**Proof.** Suppose Case 2 holds. From the continuity of the polynomial $p(x, y, z)$ on the closed and bounded parallelepiped $\Pi$, it follows that Argmin $p(\Pi)$ is a closed and bounded set. Consequently, the value $x^*$ in equality (2.6) is attained at some point $(x^*, y^*, z^*) \in$ Argmin $p(\Pi)$, and, due to bd $\Pi \cap$ Argmin $p(\Pi) = \emptyset$, we have

$$(x^*, y^*, z^*) \in \text{int}\,\Pi.$$

Moreover, since $\text{bd}\,\Pi \cap \text{Argmin}\,p(\Pi) = \emptyset$, every point in Argmin $p(\Pi)$, and hence the point $(x^*, y^*, z^*)$ as well, is a local minimum point of the polynomial $p(x, y, z)$. Then, using Statement 1, we obtain that $(x^*, y^*, z^*)$ is a local minimum point in problem (2.3) (we set $f(x, y, z) = p(x, y, z)$, $X = \Pi$, $U =$ Argmin $p(\Pi)$ in Statement 1; here, $(x^*, y^*, z^*)$ is a global minimum point in the problem $x \to \min\,;\,(x, y, z) \in U =$ Argmin $p(\Pi)$ — i.e., all conditions of Statement 1 are satisfied). Now, by Lemma 1, it follows that $x^*$ is one of the local minimum values in problem (2.2). ∎

In view of Statement 2, for the value $x^*$ satisfying (2.6), the following alternative holds:

**Case 1.** Either $x^*$ is a local minimum value in problem (2.2), which is found among the set of real roots of the univariate polynomials obtained in the alternative problems generated by applying the Variable Elimination Algorithm (VEA) to problem (2.2).

**Case 2.** Or some solution of problem (2.1) lies on bd $\Pi$ and can be found by solving one of the two-dimensional problems (analogous to the original problem (2.1), but of lower dimension):

$$p(a_1, y, z) \to \min\,;\, b_1 \le y \le b_2,\ c_1 \le z \le c_2, \tag{2.7}$$

or

$$p(a_2, y, z) \to \min\,;\, b_1 \le y \le b_2,\ c_1 \le z \le c_2, \tag{2.8}$$

or

$$p(x, b_1, z) \to \min ; a_1 \le x \le a_2,\ c_1 \le z \le c_2, \qquad (2.9)$$

or

$$p(x, b_2, z) \to \min ; a_1 \le x \le a_2,\ c_1 \le z \le c_2, \qquad (2.10)$$

or

$$p(x, y, c_1) \to \min ; a_1 \le x \le a_2,\ b_1 \le y \le b_2, \qquad (2.11)$$

or

$$p(x, y, c_2) \to \min ; a_1 \le x \le a_2,\ b_1 \le y \le b_2. \qquad (2.12)$$

To solve each of problems (2.7)–(2.12), it suffices to use the same ideas as in the treatment of the original problem (2.1).

In **Case 1**, we have the opportunity to reduce the problem dimension by one. Let $X^*$ be the finite set of real roots of the polynomials in the alternative single-variable problems with variable $x$, obtained as a result of applying the Variable Elimination Algorithm (VEA) to problem (2.2). If $X^* \cap (a_1, a_2) = \emptyset$, then, by Statement 2, the condition $\mathrm{bd}\,\Pi \cap \mathrm{Argmin}\, p(\Pi) = \emptyset$ cannot hold, and we proceed to consider **Case 2**. Otherwise, we set $X^* := X^* \cap (a_1, a_2)$. We iterate over the numbers $x^0 \in X^*$ (for example, in ascending order, starting from the minimum), and for each successive $x^0 \in X^*$, we consider a problem of the form:

$$p(x^0, y, z) \to \min \ \ (= p^*(x^0)); b_1 \le y \le b_2,\ c_1 \le z \le c_2 \qquad (2.13)$$

(i.e., a problem analogous to the original one, but with fewer variables), and we select the point $x^* = x^0 \in X^*$ with the minimal value $p^*(x^0)$, i.e.

$$p^*(x^*) = \min \{p^*(x^0) \mid x^0 \in X^*\}, x^* \in X^*.$$

Then, for any solution $(y^*, z^*)$ of problem (2.13) with $x^0 = x^*$, we obtain a solution $(x^*, y^*, z^*)$ of the original problem (2.1).

**Remark 4.** In many cases, it is also possible to find at least one global minimum point according to the algorithm from [2] (see Problem 2, pp. 112–113). Let $X^*$ be the finite set of real roots of the polynomials in the alternative single-variable problems with variable $x$, obtained as a result of applying the Variable Elimination Algorithm (VEA) to problem (2.2). If $X^* \cap (a_1, a_2) = \emptyset$, then, by Statement 2, the condition $\mathrm{bd}\ \Pi \cap \mathrm{Argmin}\ p(\Pi) = \emptyset$ cannot hold, and we proceed to consider Case 2. Otherwise, we set $X^* := X^* \cap (a_1, a_2)$. We iterate over the numbers $x^0 \in X^*$ (for example, in ascending order, starting from the minimum). For each successive $x^0 \in X^*$, we consider the sequence of alternative problems that led to the polynomial $h(x)$ in the variable $x$ (taken from one of the alternative problems), for which $x^0$ is a real root — i.e., $h(x)$ is a polynomial from an alternative problem of the form (see Remark 2):

$$x \to \min ; h(x) = 0. \qquad (2.14)$$

Let the preceding problem be

$$x \to \min; g_i(x, y) = 0,\ i = 1,2, \dots, l, \qquad (2.15)$$

where $g_i(x,y)$, $i = 1,2,\dots,l$, are polynomials — i.e., problem (2.14) is an alternative to problem (2.15) (see Remark 2). Here, problem (2.15) is preceded by problem (2.2). Then, by substituting $x = x^0$ into the system of equations

$$g_i(x^0, y) = 0, i = 1,2,\dots,l,$$

we find the possible values $y^0$ of the variable $y$. Next, by substituting $x = x^0, y = y^0$ into the system of equations in (2.2), we find the possible values $z^0$ of the variable $z$. From all possible tuples $(x^0, y^0, z^0)$, we select the point with the minimal value $p(x^0, y^0, z^0)$.

**Remark 5.** Let us analyze the obtained result. Note that the Variable Elimination Algorithm (VEA) can also be applied to problem (2.1) transformed into the form:

$$\begin{gathered} u \to \min; u - p(x,y,z) = 0, \\ (x,y,z) \in \Pi = \{(x,y,z) \in \mathbb{R}^3 | a_1 \le x \le a_2, b_1 \le y \le b_2, c_1 \le z \le c_2\}. \end{gathered} \tag{2.16}$$

In this case, we will have to perform successive elimination of three variables: $z$, $y$, and $x$. As a result, we obtain a finite number of alternative problems of the form (see Remark 2):

$$u \to \min ; h(u) = 0, \tag{2.17}$$

where $h(u)$ is a single-variable polynomial. Let $U$ be the set of all real roots of the polynomials in the alternative problems of type (2.17). Then, in the case bd $\Pi \cap$ Argmin $p(\Pi) = \emptyset$, this set will also include the value $p^* = \min\ p(\Pi)$ — the global minimum of the polynomial $p(x,y,z)$ on $\Pi$. Next, in accordance with the procedure described in Remark 4, for each $u^0 \in U$, we sequentially find the possible values of the variables $x, y, z$ by considering the corresponding sequences of alternative problems. It is quite obvious that increasing the number of variables even by one will lead to a significant complication of the considered sequence of alternative problems. Thus, a problem that can be solved using the auxiliary problem (2.2) may turn out to be practically unsolvable if the AIP is applied directly to problem (2.16).

**Remark 6.** Note that the considered case with a coordinate parallelepiped is only one of the possible cases. Suppose we need to find the global minimum in the problem:

$$\begin{gathered} p(x,y,z) \to \min ; \\ (x,y,z) \in X = \{(x,y,z) \in \mathbb{R}^3 \mid p_i(x,y,z) \le 0,\ i = 1,\dots,k\}, \end{gathered} \tag{2.18}$$

where $p(x,y,z)$, $p_1(x,y,z)$, …, $p_k(x,y,z)$ are polynomials. We assume that the polynomials $p_1(x,y,z)$, …, $p_k(x,y,z)$ have low degrees — for example, they are linear or quadratic functions. Then, similarly to the previous case, we consider the value:

$$\begin{aligned} x^* &= \min\{x \mid (x,y,z) \in X, p(x,y,z) = p^* = \min\ p(X)\} \\ &= \min\{x \mid (x,y,z) \in \text{Argmin}\ p(X)\}. \end{aligned} \tag{2.19}$$

Again, using Statement 2 and then Lemma 1, it is not difficult to show that two cases are possible:

1. $\exists (x,y,z) \in \text{Argmin}\ p(X), j \in \{1,\dots,k\}: p_j(x,y,z) = 0;$

2. $\forall(x,y,z) \in \text{Argmin } p(X)$: $p_i(x,y,z) < 0, i = 1,\dots,k,$ and then the value $x^*$ is one of the local minimum values in problem (2.2).

Note now that Case 2 is analogous to the case $\text{bd } \Pi \cap \text{Argmin } p(\Pi) = \emptyset$ for the coordinate parallelepiped $\Pi$, and it reduces to finding local minimum points in problem (2.2) (at least one of which will turn out to be a global minimum point).

In Case 1, we consider each of the auxiliary problems of the form:

$$p(x,y,z) \to \min ; p_j(x,y,z) = 0, \qquad (2.20)$$

where $j \in \{1,\dots,k\}$. Then, the global minimum point in problem (2.18) may turn out to be a local minimum point in one of the problems of type (2.20) (in the case of a single "active" constraint in problem (2.18)), that is, in the case when

$$\begin{gathered}\exists j \in \{1,\dots,k\}, \exists (x^*,y^*,z^*) \in \text{Argmin } p(X): \\ p_j(x^*,y^*,z^*) = 0, p_i(x^*,y^*,z^*) < 0, i \in \{1,\dots,k\} \backslash \{j\}.\end{gathered} \qquad (2.21)$$

Accordingly, in the case of two "active" constraints, the global minimum point in problem (2.18) may turn out to be a local minimum point in one of the problems of the form:

$$p(x,y,z) \to \min ; p_{j_1}(x,y,z) = 0,\ p_{j_2}(x,y,z) = 0,$$

where $j_1, j_2 \in \{1,\dots,k\}$, $j_1 \neq j_2$, and so on.

In the case where the polynomials $p_1(x,y,z)$, …, $p_k(x,y,z)$ have low degrees (e.g., are linear or quadratic functions), adding each new equality constraint only simplifies the work of the VEA when finding local minimum points in these problems. Even in the case of a single "active" constraint, i.e., in case (2.21), when solving a problem of the form (2.18), we can take advantage of the fact that the degrees of the variables in the polynomial $p_j(x,y,z)$ are small. This allows us to easily eliminate one of the variables from the problem

$$\begin{gathered} u \to \min, \\ u - p(x,y,z) = 0, \\ p_j(x,y,z) = 0, \end{gathered}$$

which is equivalent to problem (2.20) and to which the VEA can be applied.

**Remark 7.** As mentioned above, the case with three variables was chosen for greater clarity and simplicity of notation. The proposed method can be extended to the case with an arbitrary number of variables. However, when the number of variables is large, this possibility exists only in theory. In practice, problems with no more than 4–5 variables can be implemented.

## 3. Proof of Lemma 1

We will need some results from [6]. Let $\mathbf{x} = (x_1,\dots,x_m) \in \mathbb{R}^m$, $t \in \mathbb{R}$, and let $p_i(x_1,\dots,x_m,t)$ be polynomials in the variables $x_1,\dots,x_m,t$, for $i = 1,2,\dots,l$. Consider the one-parameter minimization problem (with parameter $t \in \mathbb{R}$):

$$\begin{cases} p_0(\mathbf{x},t) \to \min; \\ p_i(\mathbf{x},t) = 0, i = 1,2,\dots,l_1; \\ p_j(\mathbf{x},t) \le 0, j = l_1+1,\dots,l, \end{cases} \qquad (3.1)$$

(it is possible that $l_1 = 0$ and even $l = 0$). For any polynomial $p(x_1, \dots, x_m)$, denote by $\deg_{x_i} p(x_1, \dots, x_m)$ the degree of the variable $x_i$ in $p(x_1, \dots, x_m)$.

**Lemma 2.** For problem (3.1), there exists (and can even be constructively built) a polynomial $g(u,t)$ in two variables $u, t \in \mathbb{R}$ such that the degree of the variable $u$ in it is positive, and such that for any $t_0 \in \mathbb{R}$ and for any value $u_0$ of the local minimum of problem (3.1) at the fixed parameter value $t = t_0$, we have $g(u_0, t_0) = 0$.

**Proof.** Consider the problem

$$\begin{cases} u \to \min; \\ p_0(\mathbf{x},t) - u = 0; \\ p_i(\mathbf{x},t) = 0, i = 1,2,\dots,l_1; \\ p_j(\mathbf{x},t) + x_{m+j-l_1}^2 = 0, j = l_1+1,\dots,l, \end{cases} \qquad (3.2)$$

which depends on the variables $u, x_1, \dots, x_{m+l-l_1}$ and has parameter $t \in \mathbb{R}$. It is obvious that for each fixed $t \in \mathbb{R}$, the sets of local minimum values in problems (3.1) and (3.2) coincide. Applying to problem (3.2) the algorithm of successive variable elimination for $x_{m+l-l_1}, \dots, x_1$ (see [1, 2]), we obtain a finite collection of alternative problems of the form

$$u \to \min; h_i(u,t) = 0 \qquad (3.3)$$

(where $i$ is the problem index, $h_i(u,t)$ are polynomials in variables $u, t$, $i \in I$, $|I| < \infty$), such that for each $t_0 \in \mathbb{R}$ and for any number $u_0$ which, at $t = t_0$, is a local minimum value of problem (3.1) for the fixed parameter value $t = t_0$ (provided that for $t = t_0$ there exists a local minimum point in problem (3.2)), there exists an index $i_0 = i_0(u_0, t_0) \in I$ such that $u_0$ is a local minimum point of problem (3.3) for $i = i_0$, $t = t_0$, and hence $u_0$ is a root of the algebraic equation $h_i(u, t_0) = 0$ of positive degree (in the variable $u$). Then the polynomial

$$g(u,t) = \prod_{i \in I_1} h_i(u,t),$$

with $I_1 = \{i \in I \mid \deg_u h_i(u,t) > 0\}$, satisfies the requirements of the lemma if $I_1 \neq \emptyset$; and if $I_1 = \emptyset$ (in this case, problem (3.1) has no local minimum points for any fixed $t \in \mathbb{R}$), then we set, for example, $g(u,t) = u$. ∎

**Lemma 3.** Let $g(u,t)$ be a polynomial in two variables $u, t$ such that the degree of the variable $u$ in it is positive. Further, let $u_0(t)$ be a function defined and continuous on $[0, t_0]$, where $t_0 > 0$, and let $g(u_0(t), t) = 0$ hold for all $t \in [0, t_0]$. Then, for some $t_1 \in (0, t_0]$ and $r \in \mathbb{N}$, the function $u_0(t)$ can be represented as a series that converges absolutely for $t \in [0, t_0]$:

$$\sum_{i=0}^{\infty} a_i\, t^{i/r}, (a_i \in \mathbb{R},\ i = 0,1, \dots ). \qquad (3.4)$$

**Proof.** Using Theorem 21.7 from [10], we obtain that for some $t_1 \in (0, t_0]$, $r \in \mathbb{N}$, and $\delta > 0$, there exist functions $u_1(t), \dots, u_s(t)$ representable as absolutely convergent series of the form (3.4) such that for each $t \in (0, t_1]$ the following holds:

1) $u_i(t) \neq u_j(t)$ for $i \neq j$;
2) $\{u_1(t), \dots, u_s(t)\} = \{v \in \mathbb{R} \mid g(v,t) = 0,\ \mid v - u(0) \mid < \delta\}$.

Moreover, let the value $t_1$ be so small that

$$\forall t \in (0, t_1]\ \ \mid u_0(t) - u_0(0) \mid < \delta. \qquad (3.5)$$

By conditions 2) and (3.5), we have

$$\forall t \in (0, t_1]\ \ u_0(t) \in \{u_1(t), \dots, u_s(t)\}. \qquad (3.6)$$

From (3.6), condition 1), and the continuity of the functions $u_0(t), u_1(t), \dots, u_s(t)$ on $[0, t_1]$, we obtain that

$$\exists i \in \{1,2, \dots, s\}\colon u_0(t) = u_i(t),\ \forall t \in (0, t_1].$$

Indeed, suppose, for example, that $u_0(t_1) = u_1(t_1)$. We will show that

$$u_0(t) = u_1(t),\ \forall t \in (0, t_1]. \qquad (3.7)$$

Let $\bar{t} = \inf\{\tilde{t} \in (0, t_1] \mid \forall t \in [\tilde{t}, t_1]\ u_0(t) = u_1(t)\}$. If $\bar{t} = 0$, then (3.7) holds. Now suppose $\bar{t} > 0$. Then, by the continuity of the functions $u_0(t)$ and $u_1(t)$, the equality $u_0(\bar{t}) = u_1(\bar{t})$ holds, and there exists a sequence $t_n \in (0, \bar{t})$, $n = 1,2, \dots$, such that $t_n \to \bar{t}$ as $n \to \infty$ and, for example, $u_0(t_n) = u_2(t_n)$, $n = 1,2, \dots$ (here we used condition (3.6)). But then, by the continuity of the functions $u_0(t)$ and $u_2(t)$, the equality $u_0(\bar{t}) = u_2(\bar{t})$ holds, whence $u_1(\bar{t}) = u_2(\bar{t})$, which contradicts condition 1). ■

From Lemmas 2 and 3, it follows that Lemma 4 holds.

**Lemma 4.** Suppose that for each fixed $t \in (0, t_0]$, where $t_0 > 0$, problem (3.1) is solvable, and let $u(t)$ be the value of the global minimum in problem (3.1). Further, suppose that the function $u(t)$ is continuous on $(0, t_0]$ and that the limit

$$\lim_{t \to 0+} u(t) = u_0 \qquad (3.8)$$

exists. Then, for some $t_1 \in (0, t_0]$ and $r \in \mathbb{N}$, the function $u(t)$ can be represented on $(0, t_1]$ as an absolutely convergent series of the form (3.4).

**Proof.** Using Lemma 2, we obtain that for problem (3.1) there exists a polynomial $g(u,t)$ in two variables $u, t \in \mathbb{R}$ such that the degree of the variable $u$ in it is positive, and such that for any $t \in (0, t_0]$ and for the value $u(t)$ of the global minimum of problem (3.1) at the fixed parameter value $t$, the equality $g(u,t) = 0$ holds. From condition (3.8) it follows that, if we set $u(0) = u_0$, then the function $u(t)$ will be continuous on $[0, t_0]$. But then all the conditions of Lemma 3 are satisfied, from which the validity of the statement in Lemma 4 follows. ■

We will need the following strengthening of Lemma 2.

**Lemma 5.** Let $T$ be the set of numbers $t_0 \in \mathbb{R}$ such that at $t = t_0$ problem (3.1) has at least one local minimum point. Then, in the case $T \neq \emptyset$, there exists (and can even be constructively built) a polynomial $g(u,t)$ in two variables $u, t \in \mathbb{R}$ such that for any $t_0 \in T$ and for any value $u_0$ of a local minimum of problem (3.1), the quantity $u_0$ is a real root of the algebraic equation $g(u,t_0) = 0$ (with respect to the single variable $u$) of positive degree.

**Proof.** We use the collection of problems of the form (3.3) from the proof of Lemma 2. Let (as in the proof of Lemma 2) $I_1 = \{i \in I \mid \deg_u h_i(u,t) > 0\}$. Since $T \neq \emptyset$, we have $I_1 \neq \emptyset$. For each $i \in I_1$, expand the polynomial $h_i(u,t)$ in powers of the variable $u$:

$$h_i(u,t) = a_{i,s_i}(t)u^{s_i} + \cdots + a_{i,1}(t)u + a_{i,0}(t),$$

where $s_i > 0$ and $a_{i,s_i}(t) \not\equiv 0$. Set

$$q_i(t) = \mathrm{GCD}(a_{i,s_i}(t), \dots, a_{i,1}(t), a_{i,0}(t)), \tilde{a}_{i,j}(t) = \frac{a_{i,j}(t)}{q_i(t)}, j = 0,1, \dots, s_i,$$

and

$$\tilde{h}_i(u,t) = \tilde{a}_{i,s_i}(t)u^{s_i} + \cdots + \tilde{a}_{i,1}(t)u + \tilde{a}_{i,0}(t), i \in I_1.$$

Then

$$\forall i \in I_1 \quad \mathrm{GCD}\,(\tilde{a}_{i,s_i}(t), \dots, \tilde{a}_{i,1}(t), \tilde{a}_{i,0}(t)) \equiv 1, \qquad (3.9)$$

$$\forall i \in I_1,\ \forall t_0 \in T \quad \tilde{h}_i(u,t_0) \not\equiv 0. \qquad (3.10)$$

Indeed, if we suppose that for some $i_0 \in I_1$ and $t_0 \in \mathbb{R}$ we have $\tilde{h}_{i_0}(u,t_0) \equiv 0$, then $t_0$ is a real root of the polynomials

$$\tilde{a}_{i_0,s_{i_0}}(t), \dots, \tilde{a}_{i_0,1}(t), \tilde{a}_{i_0,0}(t),$$

i.e., they have a common divisor $t - t_0$, which contradicts (3.9). Denote

$$g(u,t) = \prod_{i \in I_1} \tilde{h}_i\,(u,t). \qquad (3.11)$$

We show that $g(u,t)$ satisfies the requirements of the lemma. Let $t_0 \in T$, and let $u_0$ be some value of a local minimum in problem (3.1) at $t = t_0$. Let $i_0 \in I$ be the index such that $u_0$ is a local minimum point in problem (3.3) for $i = i_0$, $t = t_0$, and hence $u_0$ is a root of the algebraic equation $h_{i_0}(u,t_0) = 0$ of positive degree. From $\deg h_{i_0}(u,t_0) > 0$ we get $\deg_u h_{i_0}(u,t) > 0$, whence $i_0 \in I_1$. In view of $h_{i_0}(u,t) = q_{i_0}(t)\tilde{h}_{i_0}(u,t)$, from $h_{i_0}(u_0,t_0) = 0$, using the fact that $\deg h_{i_0}(u,t_0) > 0$, we have $q_{i_0}(t_0) \neq 0, \tilde{h}_{i_0}(u_0,t_0) = 0,$ and deg $\tilde{h}_{i_0}(u, t_0) > 0$. Consequently, using (3.10) and (3.11), we obtain deg $g(u, t_0) > 0$ and $g(u_0, t_0) = 0$. Hence, in view of the arbitrariness of $t_0 \in T$ as well as of the value $u_0$ of a local minimum in problem (3.1) at $t = t_0$, the validity of the statement being proved follows.■

We will need the following notation. Let $X$ be a subset of $\mathbb{R}^m$, and let $y_i(x)$, $i = 1, \dots, k$, be functions defined on $X$. We introduce the sets:

$$X_0 = X, X_i = \{\mathbf{x} \in X_{i-1} \mid y_i(\mathbf{x}) = \inf\{y_i(\bar{\mathbf{x}}) \mid \bar{\mathbf{x}} \in X_{i-1}\}\}, i = 1, \dots, k.$$

Denote lexArg $[X, y_1, \dots, y_k] = X_k$ (it is possible that $X_k = \emptyset$).

We will consider power series of the form

$$q(t) = \sum_{i=1}^{\infty} a_i\, t^i, \qquad (3.12)$$

where $a_i \in \mathbb{R}$, $i = 1,2, \dots$. We assume that the series (3.12) converges absolutely for all $t \in [-t_q, t_q]$ for some $t_q > 0$. Note that for the series under consideration, we have $q(0) = 0$. The set of series of this type will be denoted by $\mathbb{R}_a\{t\}$. This set also includes the series $q(t) \equiv 0$, i.e., a series of the form (3.12) with $a_i = 0$, $i = 1,2, \dots$. Furthermore, upon substituting the series $q_1(t)$, $q_2(t)$, $q_3(t)$ into the polynomial $p(x, y, z)$, where $p(0,0,0) = 0$, the expression $p(q_1(t), q_2(t), q_3(t))$ can be transformed (after expanding the brackets and combining like terms) into a series $\bar{p}(t) = p(q_1(t), q_2(t), q_3(t))$ of the considered type, i.e.,

$$\forall q_1(t), q_2(t), q_3(t) \in \mathbb{R}_a\{t\} \;\; \bar{p}(t) = p\big(q_1(t), q_2(t), q_3(t)\big) \in \mathbb{R}_a\{t\}.$$

In addition, we will need some further concepts and notations. Let $\mathbb{R}_a^\infty\{t\}$ be the set of power series of the form

$$\sum_{i=0}^{\infty} a_i\, t^{i/r},$$

where $t \geq 0$, $a_i \in \mathbb{R}$, $i = 0,1, \dots$, $r \in \mathbb{N}$, which converge absolutely in some neighborhood of zero.

For the proof of Lemma 1, we will also need some additional statements. For $\mathbf{x} = (x_1, \dots, x_m) \in \mathbb{R}^m$, we denote

$$|\,\mathbf{x}\,| = (x_1^2 + \dots + x_m^2)^{1/2}.$$

**Lemma 6.** Let $\mathbf{x} = (x_1, \dots, x_m) \in \mathbb{R}^m$, $\mathbf{0} = (0, \dots, 0) \in \mathbb{R}^m$, and let $p_1(\mathbf{x}), \dots, p_r(\mathbf{x})$ be polynomials such that $p_1(\mathbf{0}) = 0, \dots, p_r(\mathbf{0}) = 0$. Define

$$P(t) = \{\mathbf{x} \in \mathbb{R}^m \mid p_1(\mathbf{x}) = 0, \; \dots, \; p_r(\mathbf{x}) = 0, \; |\,\mathbf{x}\,|^2 = t\}, t \geq 0,$$

and let

$$T = \{t \in \mathbb{R} \mid P(t) \neq \emptyset\}.$$

Suppose that for each $t \in T$,

$$\{\mathbf{x}^*(t)\} = \{(x_1^*(t), \dots, x_m^*(t))\} = \text{lexArg}\,[P(t), \, x_1, \; \dots, \, x_m].$$

Then $T \neq \emptyset$, and for each $i \in \{1,2, \dots, m\}$ there exists a polynomial $g_i(x_i, t)$ in the variables $x_i$ and $t$ such that for any $t_0 \in T$ the value $x_i^*(t_0)$ is a root of the algebraic equation $g_i(x_i, t_0) = 0$ (in the variable $x_i$) of positive degree.

**Proof.** Note that $0 \in P(0)$, whence $0 \in T$ and thus $T \neq \emptyset$. Consider the parametric global minimization problem (with parameter $t \in T \neq \emptyset$):

$$x_1 \to \min\ (= x_1^*(t)); p_1(\mathbf{x}) = 0,\ \dots,\ p_r(\mathbf{x}) = 0, |\ \mathbf{x}\ |^2 = t. \qquad (3.13)$$

By Lemma 5, for problem (3.13) there exists a polynomial $g_1(x_1, t)$ in the variables $x_1$ and $t$ such that for any $t_0 \in T$ the value $x_1^*(t_0)$ is a root of the algebraic equation (in the variable $x_1$) of positive degree

$$g_1(x_1, t_0) = 0. \qquad (3.14)$$

Next, consider the problem

$$x_2 \to \min; p_1(\mathbf{x}) = 0, \dots,\ p_r(\mathbf{x}) = 0, |\ \mathbf{x}\ |^2 = t,\ g_1(x_1, t) = 0, \qquad (3.15)$$

as well as the global minimization problem

$$\begin{gathered} x_2 \to \min\ (= x_2^*(t)); \\ p_1(x_1^*(t), x_2, \dots, x_m) = 0, \dots,\ p_r(x_1^*(t), x_2, \dots, x_m) = 0; \\ |\ (x_1^*(t), x_2, \dots, x_m)\ |^2 = t, \end{gathered} \qquad (3.16)$$

where $t \in T$ is the parameter in problems (3.15) and (3.16). We will show that if $(x_2^*(t_0), \check{x}_3(t_0), \dots, \check{x}_m(t_0))$ is, for some $t_0 \in T$, a global minimum point in problem (3.16) at $t = t_0$, then

$$(x_1^*(t_0), x_2^*(t_0), \check{x}_3(t_0), \dots, \check{x}_m(t_0)) \qquad (3.17)$$

is a local minimum point in problem (3.15) at $t = t_0$. We use the fact that $x_1^*(t_0)$ is a root of the algebraic equation $g_1(x_1, t_0) = 0$ (in the variable $x_1$) of positive degree. Let $\varepsilon > 0$ be such that in the $\varepsilon$-neighborhood of the point $x_1^*(t_0)$ there are no other roots of equation (3.14). Then

$$\forall v \in \mathbb{R} \quad [\{(|v - x_1^*(t_0)| \leq \varepsilon) \wedge (g_1(v, t_0) = 0)\} \Rightarrow v = x_1^*(t_0)]. \quad (3.18)$$

If we suppose that point (3.17) is not a local minimum point in problem (3.15) at $t = t_0$, then there exists a point $\tilde{\mathbf{x}} = (\tilde{x}_1, \tilde{x}_2, \dots, \tilde{x}_m)$ lying in the $\varepsilon$-neighborhood of point (3.17), such that

$$p_1(\tilde{\mathbf{x}}) = 0,\ \dots,\ p_r(\tilde{\mathbf{x}}) = 0,\ |\ \tilde{\mathbf{x}}\ |^2 = t_0,\ g_1(\tilde{x}_1, t_0) = 0,$$

and

$$\tilde{x}_2 < x_2^*(t_0). \qquad (3.19)$$

Then from (3.18) we obtain that $\tilde{x}_1 = x_1^*(t_0)$, whence

$$\begin{gathered} p_1(x_1^*(t_0), \tilde{x}_2, \dots, \tilde{x}_m) = 0,\ \dots,\ p_r(x_1^*(t_0), \tilde{x}_2, \dots, \tilde{x}_m) = 0; \\ |\ (x_1^*(t_0), \tilde{x}_2, \dots, \tilde{x}_m)\ |^2 = t_0. \end{gathered} \qquad (3.20)$$

Note now that conditions (3.19) and (3.20) contradict the fact that $(x_2^*(t_0), \check{x}_3(t_0), \dots, \check{x}_m(t_0))$ is a global minimum point in problem (3.16) at $t = t_0$.

Using the result proved above, by Lemma 5 (we apply Lemma 5 to problem (3.15)) we obtain that there exists a polynomial $g_2(x_2, t)$ in the variables $x_2$ and $t$ such that for any $t_0 \in T$ the value $x_2^*(t_0)$ is a root of the algebraic equation $g_2(x_2, t_0) = 0$ (in the variable $x_2$) of positive degree.

Next, similarly to problems (3.15) and (3.16), we consider the problem

$$x_3 \to \min; p_1(\mathbf{x}) = 0, \dots, p_r(\mathbf{x}) = 0, |\mathbf{x}|^2 = t, g_1(x_1, t) = 0, g_2(x_2, t) = 0, (3.21)$$

as well as the problem

$$\begin{gathered} x_3 \to \min \ \ (= x_3^*(t)); \\ p_1(x_1^*(t), x_2^*(t), x_3, \dots, x_m) = 0, \ \dots ; \\ p_r(x_1^*(t), x_2^*(t), x_3, \dots, x_m) = 0; \\ |\ (x_1^*(t), x_2^*(t), x_3, \dots, x_m)\ |^2 = t, \end{gathered} \tag{3.22}$$

where $t \in T$ is the parameter in problems (3.21) and (3.22). Then, arguing analogously to the analysis of problems (3.15) and (3.16), we obtain that if $(x_3^*(t_0), \check{x}_4(t_0), \dots, \check{x}_m(t_0))$ is, for some $t_0 \in T$, a global minimum point in problem (3.22) at $t = t_0$, then

$$(x_1^*(t_0), x_2^*(t_0), x_3^*(t_0), \check{x}_4(t_0), \dots, \check{x}_m(t_0))$$

is a local minimum point in problem (3.21) at $t = t_0$. But then, similarly to the previous reasoning, using Lemma 5 (we apply Lemma 5 to problem (3.21)), we obtain that there exists a polynomial $g_3(x_3, t)$ in the variables $x_3$ and $t$such that for any $t_0 \in T$ the value $x_3^*(t_0)$ is a root of the algebraic equation $g_3(x_3, t_0) = 0$ (in the variable $x_3$) of positive degree.

Continuing our reasoning, at the subsequent stages we will successively obtain the required polynomials $g_4(x_4, t), \dots, g_m(x_m, t)$. ■

**Theorem 1.** Let $\mathbf{x} = (x_1, \dots, x_m) \in \mathbb{R}^m$, and let $\mathbf{0} = (0, \dots, 0) \in \mathbb{R}^m$ be a local minimum point of the polynomial $p(\mathbf{x})$, with $p(\mathbf{0}) = 0$. Suppose that $\mathbf{0}$ is simultaneously a local minimum point in the problem

$$x_1 \to \min ; p(\mathbf{x}) = 0.$$

Then $\mathbf{0}$ is a local minimum point in the problem

$$x_1 \to \min ; \frac{\partial p(\mathbf{x})}{\partial x_i} = 0, \ i = 1, \dots, m. \tag{3.23}$$

**Proof.** Suppose that $\mathbf{0}$ is not a local minimum point in problem (3.23). Consider the set depending on the parameter $t \geq 0$:

$$P(t) = \{\mathbf{x} \in \mathbb{R}^m \ | \ \frac{\partial p(\mathbf{x})}{\partial x_i} = 0, \ i = 1, \dots, m, \ |\mathbf{x}|^2 = t\}.$$

Let

$$T = \{t \in \mathbb{R} \mid P(t) \neq \emptyset\},$$

and for each $t \in T$, let

$$\{\mathbf{x}^*(t)\} = \{(x_1^*(t), \dots, x_m^*(t))\} = \operatorname{lexArg}\,[P(t), x_1, \dots, x_m].$$

Note that $0 \in T$, whence $T \neq \emptyset$. Since for each $t \in T$ the set $P(t)$ is nonempty, bounded, and closed, it follows that for each $t \in T$ we have

$$\operatorname{lexArg}\,[P(t), \ x_1, \dots, x_m] \neq \emptyset.$$

Because $\mathbf{0}$ is not a local minimum point in problem (3.23), there exists a sequence of numbers

$$t_n \in T, t_n > 0, \lim_{n\to\infty} t_n = 0, P(t_n) \neq \emptyset, x_1^*(t_n) < 0, n = 1,2, \dots . \quad (3.24)$$

By Lemma 6, there exist polynomials $g_i(x_i, t)$ in the variables $x_i$ and $t$ such that $\deg_{x_i} g_i(x_i, t) > 0$ for $i = 1, \dots, m$, and

$$\forall t \in T \quad g_i(x_i^*(t), t) = 0, i = 1, \dots, m.$$

Using Theorem 21.7 from [10], we obtain that there exist series

$$q_{ij}(t) \in \mathbb{R}_a^\infty\{t\}, j = 1,2, \dots, s_i \in \mathbb{N}, \ i = 1, \dots, m, \quad (3.25)$$

such that for some $t_1 > 0$ and $\delta > 0$ the series (3.25) converge absolutely on $[0, t_1]$, and

$$\forall t \in [0, t_1] \ \{x \in \mathbb{R} | \, |x| \leq \delta, g_i(x, t) = 0\} =$$
$$= \{v \in \mathbb{R} | v = q_{ij}(t), j = 1,2, \dots, s_i\}, i = 1, \dots, m.$$

But then, for the sake of notational simplicity and to avoid passing to a subsequence, we may assume that

$$x_i^*(t_n) = q_{i1}(t_n), i = 1, \dots, m, \ n = 1,2, \dots .$$

At the same time, since $|\mathbf{x}^*(t_n)|^2 = t_n$, we have $x_i^*(t_n) \to 0$ as $n \to \infty$ for $i = 1, \dots, m$, whence

$$q_{i1}(0) = 0, i = 1, \dots, m. \quad (3.26)$$

In view of (3.25) and (3.26), for some $k \in \mathbb{N}$ we have

$$\bar{q}_i(\tau) = q_{i1}(\tau^k) \in \mathbb{R}_a\{\tau\}, i = 1, \dots, m,$$

and, since $\mathbf{x}^*(t_n) \in P(t_n)$, taking (3.24) into account, we obtain:

$$\frac{\partial p(\bar{q}_1(\tau_n), \dots, \bar{q}_m(\tau_n))}{\partial x_i} = 0, i = 1, \dots, m,$$

$$\left|\left(\bar{q}_1(\tau_n), \dots, \bar{q}_m(\tau_n)\right)\right|^2 = t_n = \tau_n^k > 0, \bar{q}_1(\tau_n) = q_{11}(t_n) = x_1^*(t_n) < 0, \quad (3.27)$$

where $\tau_n = t_n^{1/k}$, $n = 1,2, \dots$. But then (see Remark 9)

$$\frac{\partial p(\bar{q}_1(\tau), \dots, \bar{q}_m(\tau))}{\partial x_i} \equiv 0, i = 1, \dots, m,$$

whence

$$\frac{d}{d\tau} p(\bar{q}_1(\tau), \dots, \bar{q}_m(\tau)) = \frac{\partial p(\bar{q}_1(\tau), \dots, \bar{q}_m(\tau))}{\partial x_1} \cdot \frac{d\bar{q}_1(\tau)}{d\tau} + \cdots$$
$$+ \frac{\partial p(\bar{q}_1(\tau), \dots, \bar{q}_m(\tau))}{\partial x_m} \cdot \frac{d\bar{q}_m(\tau)}{d\tau} \equiv 0.$$

And, since (see (3.26)) we have $p(\bar{q}_1(0), \dots, \bar{q}_m(0)) = p(\mathbf{0}) = 0$, it follows that

$$p(\bar{q}_1(\tau), \dots, \bar{q}_m(\tau)) \equiv 0,$$

whence

$$p(\bar{q}_1(\tau_n), \dots, \bar{q}_m(\tau_n)) = 0, n = 1,2, \dots .$$

Taking (3.27) into account, this contradicts the fact that $\mathbf{0}$ is a local minimum point in the problem

$$x_1 \to \min ; p(\mathbf{x}) = 0. \blacksquare$$

**Remark 8.** The statement of Theorem 1 can be obviously extended to the case where, instead of the variable $x_1$, any of the variables $x_i$, $i = 1, \dots, m$, is chosen in its conditions. Moreover, instead of $x_i \to \min$, one can use $-x_i \to \min$, where $i \in \{1, \dots, m\}$.

**Remark 9.** Here we use the following simple assertion. If $q(t) \in \mathbb{R}_a\{t\}$, $t_n > 0$, the series $q(t_n)$ converges absolutely, and $q(t_n) = 0$ for $n = 1,2, \dots$, with $\lim_{n\to\infty} t_n = 0$, then $q(t) \equiv 0$.

**Remark 10.** Lemma 1 is an obvious corollary of Theorem 1.

We also present one of the corollaries of Theorem 1. We will need the following simple statement.

**Statement 3.** Let $X$ be a nonempty subset of $\mathbb{R}^m$, and let $\mathbf{x}^* \in X$. Then, in order for $\mathbf{x}^*$ to be an isolated point of the set $X$ (i.e., there exists $\delta > 0$ such that $B(\mathbf{x}^*, \delta) \cap X = \{\mathbf{x}^*\}$, where $B(\mathbf{x}^*, \delta) = \{x \in \mathbb{R}^m \mid | \mathbf{x} - \mathbf{x}^* | \le \delta\}$), it is necessary and sufficient that $\mathbf{x}^*$ be a local minimum point in the problem

$$x_i \to \min ; \mathbf{x} \in X,$$

as well as in the problem (see Remark 8)

$$-x_i \to \min ; \mathbf{x} \in X,$$

for all $i = 1,2, \dots, m$.

As a corollary of Theorem 1 and Statement 3, we have

**Theorem 2.** Let $\mathbf{0} = (0, \dots, 0) \in \mathbb{R}^m$ be an isolated local minimum point of the polynomial $p(x)$, where $x \in \mathbb{R}^m$. Then $\mathbf{0}$ is an isolated point of the set

$$\left\{\mathbf{x} \in \mathbb{R}^m \middle| \frac{\partial p(\mathbf{x})}{\partial x_i} = 0, i = 1, \dots, m\right\}.$$

## References

1. Nefedov, V. N. (1987). Polynomial optimization problems. *USSR Computational Mathematics and Mathematical Physics*, *27*(3), 13–21. https://doi.org/10.1016/0041-5553(87)90074-7
2. Nefedov, V. N., & Zharkikh, A. V. (2020). Algorithmization and software implementation of the variable elimination method in polynomial optimization problems. *Modeling and Data Analysis*, *10*(1), 110–128. https://doi.org/10.17759/mda.2020100107
3. Lasserre, J. B. (2001). Global optimization with polynomials and the problem of moments. *SIAM Journal on Optimization*, *11*(3), 796–817.
4. Laurent, M. (2009). Sums of squares, moment matrices and optimization over polynomials. In *Emerging Applications of Algebraic Geometry* (*IMA Volumes in Mathematics and its Applications*, Vol. 149, pp. 157–270). Springer.
5. Hansen, E. R. (1992). *Global optimization using interval analysis*. New York: M. Dekker. (Monographs and Textbooks in Pure and Applied Mathematics, Vol. 165).

6. Kearfott, R. B. (1992). An interval branch and bound algorithm for bound constrained optimization problems. *Journal of Global Optimization*, *2*(3), 259–280.
7. Lebbah, Y. (2009). ICOS: a branch and bound based solver for rigorous global optimization. *Optimization Methods & Software*, *24*(4–5), 609–626.
8. [Author(s) not specified]. (2016). Interval branch-and-bound algorithms for optimization and constraint satisfaction: a survey and prospects. *Journal of Global Optimization*, *65*, 837–866.
9. Nefedov, V. N. (1990). Estimation of the error in convex polynomial optimization problems. *U.S.S.R. Computational Mathematics and Mathematical Physics*, *30*(1), 147–158. https://doi.org/10.1016/0041-5553(90)90024-M
10. Krasnosel'skii, M. A., Vainikko, G. M., Zabreiko, P. P., et al. (1972). *Approximate solution of operator equations*. Leiden, The Netherlands: Noordhoff.